# Proof of Fermat's Last Theorem by Algebra Identities and Linear Algebra

Javad Babaee Ragani*

*Young Researchers and Elite Club, Qaemshahr Branch, Islamic Azad University, Qaemshahr, Iran
*Department of Civil Engineering, Arak Branch, Islamic Azad University, Arak, Iran
Email: j.babaee@gaemiau.ac.ir

**Abstract**

The main aim of the present paper is to represent an exact and simple proof for Fermat's Last Theorem by using properties of the algebra identities and linear algebra.

**1. Introduction and preliminaries**

In number theory, Fermat's Last Theorem states that: no three positive integers $a$, $b$ and $c$ can satisfy the equation $a^x + b^x = c^x$ for any integer value of $x$ greater than two. This theorem was first conjectured by Pierre de Fermat in 1637 and the first successful proof was released in 1994 by Professor Andrew Wiles and finally published in 1995 after 358 years of efforts by mathematicians but the proof was well over 100 pages long and complex utilizing the most modern 20$^{th}$ century analytical although Fermat claimed that: "I've found a remarkable proof of this fact, but there is not enough space in the margin [of the book] to write it". Therefore a question has preoccupied the minds of mathematicians always and that is, is there a 'simple' proof of the Theorem? The author believes that based on the Lemmas and proving presented in this article, maybe he achieved to a simple proving actually.

**Lemma 1.1.** It is suffices to prove Fermat's Last Theorem for 4 and for every odd prime $p \geq 3$.

It is well known that if the Last Theorem can be proved for $n = 4$, then it is also proven for all multiples of $n = 4$, because all of the remaining numbers can be reduced to a multiple of the prime numbers, it is therefore only necessary to prove Fermat's Last theorem for all the primes.

**Lemma 1.2.** In equation $a^{2k+1} + b^{2k+1} = c^{2k+1}$, the expressions $(a+b), (c-b)$, and $(c-a)$ are coprime.

**Proof.** In equation $a^{2k+1} + b^{2k+1} = c^{2k+1}$ numbers $a$, $b$, and $c$ are relatively prime in pairs and because that equation can be written in the form of three relations (1) to (3), so it can be concluded that the terms $(a+b)$, $(c-a)$, and $(c-b)$ must be relatively prime in pairs.

$$a^{2k+1} + b^{2k+1} = (a+b)\left(a^{2k} - ab^{2k-1} + a^2b^{2k-2} - \ldots \pm (ab)^k\right) = c^{2k+1} \qquad (1)$$

$$c^{2k+1} - a^{2k+1} = (c-a)\left(a^{2k} + ac^{2k-1} + a^2c^{2k-2} - \ldots + (ac)^k\right) = b^{2k+1} \qquad (2)$$

$$c^{2k+1} - b^{2k+1} = (c-b)\left(b^{2k} + bc^{2k-1} + b^2c^{2k-2} - \ldots + (bc)^k\right) = a^{2k+1} \qquad (3)$$

**Theorem 1.3.** Consider the equation $a^{2k+1} + b^{2k+1} = c^{2k+1}$ which in that $gcd(a, b, c) = 1$ and $2k+1$ is prime number and let $a = a_1 a_2 a'$, $b = b_1 b_2 b'$, and $c = c_1 c_2 c'$, if set up above equation then will be set up three following couple-relations simultaneously:

$$a+b = c_1^{2k+1} * c'^{2k}, \quad a^{2k} + b^{2k} - ab^{2k-1} - ba^{2k-1} \ldots \pm (ab)^k = c_2^{2k+1} * c'$$

$$c-a = b_1^{2k+1} * b'^{2k}, \quad c^{2k} + a^{2k} + ac^{2k-1} + ca^{2k-1} \ldots + (ac)^k = b_2^{2k+1} * b'$$

$$c-b = a_1^{2k+1} * a'^{2k}, \quad c^{2k} + b^{2k} + cb^{2k-1} + bc^{2k-1} \ldots + (bc)^k = a_2^{2k+1} * a'$$

That we have two case in above relations, case I: $a' = b' = c' = 1$ $(2k+1 \nmid abc)$ and case II: exactly one of the numbers $a'$, $b'$ or $c'$ is equal to $2k+1$ $(2k+1 | abc)$.

**Proof.** Consider the equation $a^{2k+1} + b^{2k+1} = c^{2k+1}$ which in that positive numbers a, b, and c are relatively prime in pair and $2k+1$ is prime number. Now let $c = c_1 * c_2 * c_1' * c_2' * \ldots * c_n'$, therefore the expression $(a^{2k+1} + b^{2k+1})$ based on the properties of algebra identities can be written as follows:

$$a^{2k+1} + b^{2k+1} = (a+b)(a^{2k} + b^{2k} + \ldots \pm (ab)^k) = c^{2k+1} = c_1^{2k+1} * c_2^{2k+1} * c_1'^{2k+1} * c_2'^{2k+1} * \ldots * c_n'^{2k+1} \quad (4)$$

The equation (4) can be written in the two following forms:

$$a + b = c_1^{2k+1} * c_1'^{m_1} * c_2'^{m_2} * \ldots * c_n'^{m_n} \quad (5)$$

$$a^{2k} + b^{2k} - ab^{2k-1} - ba^{2k-1} \ldots \pm (ab)^k = c_2^{2k+1} * c_1'^{2k+1-m_1} * c_2'^{2k+1-m_2} * \ldots * c_n'^{2k+1-m_n} \quad (6)$$

That in the above relations, exponent $c_1$ is multiples of $c$ that term $(a+b)$ just dividable by them ($c_1 | (a+b)$ but $c_1 \nmid [a^{2k} + b^{2k} - ab^{2k-1} - ba^{2k-1} \ldots \pm (ab)^k]$) and exponent $c_2$ is multiples of $c$ that term $[a^{2k} + b^{2k} - ab^{2k-1} - ba^{2k-1} \ldots \pm (ab)^k]$ just dividable by them ($c_2 | [a^{2k} + b^{2k} - ab^{2k-1} - ba^{2k-1} \ldots \pm (ab)^k]$ but $c_2 \nmid (a+b)$) and exponents $c_1'$ to $c_n'$ are multiples of $c$ that both term $(a+b)$ and $[a^{2k} + b^{2k} - ab^{2k-1} - ba^{2k-1} \ldots \pm (ab)^k]$ are dividable by them and in other words are the common factors of mentioned terms ($gcd(a+b, a^{2k} + b^{2k} - ab^{2k-1} - ba^{2k-1} \ldots \pm (ab)^k) = c_1' * c_2' * \ldots * c_n'$), therefore $c_1$ and $c_2$ are relatively prime. Relation (6) can be expanded by the use of $(a+b)$ of relation (5) and after summarizing can be written in the form of relation (9):

$$a^{2k} + b^{2k} - ab^{2k-1} - ba^{2k-1} + \ldots \pm (ab)^k = (a+b)^{2k} + q_1(a+b)^{2k-2} + q_2(a+b)^{2k-4} + \ldots$$
$$+ q_{k-1}(a+b)^2 \pm (2k+1)(ab)^k = (a+b)^2[(a+b)^{2k-2} + q_1(a+b)^{2k-4} + q_2(a+b)^{2k-6} + \quad (7)$$
$$\ldots + q_{k-1}] \pm (2k+1)(ab)^k = c_2^{2k+1} * c_1'^{2k+1-m_1} * c_2'^{2k+1-m_2} * \ldots * c_n'^{2k+1-m_n} \Rightarrow$$

$$\pm (2k+1)(ab)^k = (c_2^{2k+1} * c_1'^{2k+1-m_1} * c_2'^{2k+1-m_2} * \ldots * c_n'^{2k+1-m_n}) - \{\overbrace{(c_1^{2k+1} * c_1'^{m_1} * c_2'^{m_2} * \ldots * c_n'^{m_n})^2}^{(a+b)^2} *$$
$$\underbrace{[(a+b)^{2k-2} + q_1(a+b)^{2k-4} + q_2(a+b)^{2k-6} + \ldots + q_{k-1}]}_{M}\} \Rightarrow \quad (8)$$

$$\pm (2k+1)(ab)^k = c_2^{2k+1} * c_1'^{2k+1-m_1} * \ldots * c_n'^{2k+1-m_n} - c_1^{2(2k+1)} * c_1'^{2m_1} * c_2'^{2m_2} \ldots * c_n'^{2m_n} * M \quad (9)$$

That in above relations, M is integer and amounts $q_1$ to $q_{k-1}$ are definite integer numbers (see the example 1.4). In both right term of the relation (9) there are the common factors $c_1'$ to $c_n'$ with a different that in one of the expressions have the powers $(2k+1-m_1)$ to $(2k+1-m_n)$ respectively and in other expressions have the powers $2m_1$ to $2m_n$ respectively. Since this factors $(c_i')$ themselves are a multiple of $c$, so in comparison with the term $(ab)^k$ in the left side of the relation (9) are prime $\left(gcd\left(ab, c_1' * c_2' * ... * c_n'\right) = 1\right)$ and as the number $2K+1$ considered to be a prime number, therefore this number have just a factor (itself only that is $2K+1$) and thus among the common factors $c_1'$ to $c_n'$ just one of those numbers have an definite amount and the rest of them must be equal to 1 necessarily. Suppose that $c_1'$ has a definite amount and therefore always $c_2'$ to $c_n'$ are equal to 1. Thus the equation (9) can be written as follows:

$$\pm(2k+1)(ab)^k = c_2^{2k+1} * c_1'^{2k+1-m_1} - c_1^{2(2k+1)} * c_1'^{2m_1} * M \qquad (10)$$

Now it is tired to declare that $c_1'$ is equal to 1 or equal to $2K+1$. In relation (10) for the power $c_1'$, two states can be considered: or state I: $(2k+1-m_1) > 2m_1$ and or state II: $(2k+1-m_1) < 2m_1$. In the first state that is $(2k+1-m_1) > 2m_1$ the relation (10) can be written according to the following relation:

$$\pm(2k+1)(ab)^k = c_1'^{2m_1} * [(c_2^{2k+1} * c_1'^{2k+1-3m_1}) + (c_1^{2(2k+1)} * M)] \qquad (11)$$

In relation (11) as it referred before $gcd\left(ab, c_1'\right) = 1$, therefore $c_1'^{2m_1}$ can be equal to 1 or $2k+1$. Because $2k+1$ is prime thus it's power always is equals to 1. To make an equality of $(2k+1)^1 = c_1'^{2m_1}$ there must be $c_1' = 2k+1$ and $2m_1 = 1$ that is a contradiction, so in the first case there must be $c_1' = 1$. Now the second state for the powers is supposed, it's means $(2k+1-m_1) < m_1$. In this case the relation (11) can be written as follows:

$$\pm(2k+1)(ab)^k = c_1'^{2k+1-m_1} * [c_2^{2k+1} + (c_1'^{-2k-1+3m_1} * c_1^{2(2k+1)} * M)] \qquad (12)$$

In relation (12) based on above concepts there must be: $2k+1 = c_1'^{2k+1-m_1}$ that to set up this equality we must have $c_1' = 2k+1$ and $m_1 = 2k$ (because $2k+1-m_1 = 1$). So to study the relation (10) by assuming that $2k+1$ is a prime number it was indicated that firstly the amounts $c_2'$ to $c_n'$ must be equal to 1 and secondly by considering the relation (10) and the two states which can be for the power of $c_1'$ it was proved that $c_1'$ is equal to 1 or $2k+1$ that in the second state ($c_1' = 2k+1$) the amount of power $m_1$ must be equal to $2k$ and so $c_1$ and $c_2$ are coprime. Thus the couple-relations (5)-(6) can be written in the final form of the following couple-relations:

$$a + b = c_1^{2k+1} * c'^{2k} \qquad (13)$$
$$a^{2k} + b^{2k} - ab^{2k-1} - ba^{2k-1} ... \pm (ab)^k = c_2^{2k+1} * c' \qquad (14)$$

That as it mentioned before, in the above couple-relations $c'$ that is greatest common divisor (gcd) terms $(a+b)$ and $(a^{2k} + b^{2k} - ab^{2k-1} - ba^{2k-1} ... \pm (ab)^k$ is equal to 1 or equal to $2k+1$. Now if the equation $a^{2k+1} + b^{2k+1} = c^{2k+1}$ will be written in the form of equations $c^{2k+1} + b^{2k+1} = -a^{2k+1}$ and $c^{2k+1} + a^{2k+1} = -b^{2k+1}$ based on the reasoning mentioned above and conducting the same procedure it can be indicated that above equations can be written according to the couple-relations (15)-(16) and (17)-(18) respectively:

$$c - a = b_1^{2k+1} * b'^{2k} \tag{15}$$

$$c^{2k} + a^{2k} + ac^{2k-1} + ca^{2k-1} \ldots + (ac)^k = b_2^{2k+1} * b' \tag{16}$$

$$c - b = a_1^{2k+1} * a'^{2k} \tag{17}$$

$$c^{2k} + b^{2k} + bc^{2k-1} + cb^{2k-1} \ldots + (bc)^k = a_2^{2k+1} * a' \tag{18}$$

In each of the above couple-relations the greatest common dividers $a'$ and $b'$ are either equal to 1 or equal to $2k+1$. Based on Lemma 1.2 as the terms $(a+b)$, $(c-b)$, and $(c-a)$ are relatively prime in pairs, so exactly one of the common factors $a'$, $b'$, or $c'$ can be equal to $2k+1$ and the other two numbers are equal to 1 necessarily, in other words two states can be considered for the above couples-relations: case I: ($2k+1 \nmid abc$) or none of the above couple-relations has no common factor that is $a' = b' = c' = 1$ and or case II: ($2k+1 \mid abc$) or exactly one of the above couple-relations has a common factor that is or $a' = 2k+1$ (and $b' = c' = 1$) or $b' = 2k+1$ (and $a' = c' = 1$) and or $c' = 2k+1$ (and $a' = b' = 1$).

**Example 1.4.** Example of Theorem 1.3 for case $2K + 1 = 7$

$$c^7 - a^7 = b^7 = (b_1 * b_2 * b'_1 * b'_2 * \ldots * b'_n)^7$$

That $b_1$, $b_2$, and $b'_1$ to $b'_n$ are relatively prime in pair and exponent $b_1$ just counts term $(c-a)$, (but $\gcd(a^6 + c^6 + a c^5 + c a^5 + a^2 c^4 + c^4 a^2 + a^3 c^3, b_1) = 1$) and exponent $b_2$ just counts term $(a^6 + c^6 + a c^5 + c a^5 + a^2 c^4 + c^4 a^2 + a^3 c^3)$, (but $\gcd(c-a, b_2) = 1$) and $b'_1$ to $b'_n$ are the common factors of both term.

$$c - a = b_1^7 * b'^{m_1}_1 * b'^{m_2}_2 * \ldots * b'^{m_n}_n \quad (m_1 \text{ to } m_n \text{ are indefinite powers})$$

$$a^6 + c^6 + ac^5 + ca^5 + a^2c^4 + c^4a^2 + (ac)^3 = b_2^7 * b'^{7-m_1}_1 * b'^{7-m_2}_2 * \ldots * b'^{7-m_n}_n =$$

$$(c-a)^6 + 7ac(c-a)^4 + 14(ac)^2(c-a)^2 + 7(ac)^3 \Rightarrow$$

$$7(ac)^3 = b_2^7 * b'^{7-m_1}_1 * b'^{7-m_2}_2 * \ldots * b'^{7-m_n}_n - (c-a)^2 \left[ (c-a)^4 + 7ab(c-a)^2 + 14(ac)^2 \right] =$$

$$b_2^7 * b'^{7-m_1}_1 * \ldots * b'^{7-m_n}_n - (b_1^{14} * b'^{2m_1}_1 * b'^{2m_2}_2 * \ldots * b'^{2m_n}_n) \left[ (c-a)^4 + 7ab(c-a)^2 + 14(ac)^2 \right] \Rightarrow$$

$$7(ac)^3 = b_2^7 * b'^{7-m_1}_1 * b'^{7-m_2}_2 * \ldots * b'^{7-m_n}_n - (b_1^{14} * b'^{2m_1}_1 * b'^{2m_2}_2 * \ldots * b'^{2m_n}_n) * M \Rightarrow$$

Because exponents $b'_i$ and $ac$ are coprime and 7 is prime number, thus either $b'_1 = 1$ and or $b'_1 = 2k+1$ and so $b'_2$ to $b'_n$ are equals to 1 necessarily. $\Rightarrow$

$$7(ac)^3 = b_2^7 * b'^{7-m_1}_1 - (b_1^{14} * b'^{2m_1}_1) * M \Rightarrow$$

Either $b'^{2m_1}_1 = 1$ or $7 \rightarrow$ If so just $b'_1 = 1$ and or $b'^{7-m_1}_1 = 1$ or $7 \rightarrow b'_1 = 1$ or $b'_1 = 7$ (In this case: $7 - m_1 = 1$ or $m_1 = 6$)

$$c^7 - a^7 = b^7 \rightarrow c - a = b_1^7 * b'^6_1 \text{ and } a^6 + c^6 + ac^5 + ca^5 + a^2c^4 + c^4a^2 + (ac)^3 = b_2^7 * b'_1$$

That $b_1'$ is equal 1 or 7. If in relation $a^7 + b^7 = c^7$ it is supposed that $a = a_1 a_2 a'$ and $b = b_1 b_2 b'$ then according to proving process mentioned above and in a similar way it can be declared that the couple-relations $a + b = c_1^7 . c_1'^6$,

$a^6 + b^6 + ba^5 + ab^5 + b^2 a^4 + a^4 b^2 + (ba)^3 = c_2^7 . c'$ and so couple-relations $c - b = a_1^7 . a_1'^6$,

$b^6 + c^6 + bc^5 + cb^5 + b^2 c^4 + c^4 b^2 + (bc)^3 = a_2^7 . a'$ must be set up.

Therefore it is generally represented that in relation $a^{2k+1} + b^{2k+1} = c^{2k+1}$ when $2k+1$ is a prime number, can be written in the form of there couple-relations (13)-(14), (15)-(16), and (17)-(18) and in other words if that relation is set up, the above three couple-relations must be set up simultaneously that necessarily two states can exist for the common numbers $a'$, $b'$, and $c'$, case I: ($2k+1 \nmid abc$) or $a' = b' = c' = 1$ and or case II: ($2k+1 | abc$) or exactly one of the numbers $a'$ or $b'$ or $c'$ is equal to $2k+1$ (and the other two numbers should be equal to one necessarily).

**Lemma 1.5.** In equation $a^{2k+1} + b^{2k+1} = c^{2k+1}$, amount expression $(a+b-c)$ is always positive

**Proof.** $(a+b)^{2k+1} \geq a^{2k+1} + b^{2k+1} \geq c^{2k+1} \rightarrow a+b \geq c \rightarrow a+b-c \geq 0$

**Lemma 1.6.** In equation $a^{2k+1} + b^{2k+1} = c^{2k+1}$ in the lieu of each positive $a$, $b$, and $c$ the relations $a+b-c = a_1 a'(a_2 - a_1^{2k} a'^{2k-1}) = b_1 b'(b_2 - b_1^{2k} b'^{2k-1}) = c_1 c'(c_1^{2k} c'^{2k-1} - c_2) \geq 0$ is always set up

**Proof.** According to Theorem 1.4 if the relation $a^{2k+1} + b^{2k+1} = c^{2k+1}$ is set up assuming $a = a_1 a_2 a'$ and $b = b_1 b_2 b'$ and $c = c_1 c_2 c'$, then the relations (13), (15), and (17) will be set up simultaneously and in this case the amount of the expression $(a+b-c)$ can be written as follows:

$(a+b) - c = c_1^{2k+1} c'^{2k} - c_1 c_2 c' = c_1 c'(c_1^{2k} c'^{2k-1} - c_2) \geq 0 \rightarrow c_1 c' > 0$ thus $c_1^{2k} c'^{2k-1} - c_2 \geq 0$ (19)

$a + (b-c) = a_1 a_2 a' - a_1^{2k+1} a'^{2k} = a_1 a'(a_2 - a_1^{2k} a'^{2k-1}) \geq 0 \rightarrow a_1 a' > 0$ thus $a_2 - a_1^{2k} a'^{2k-1} \geq 0$ (20)

$b + (a-c) = b_1 b_2 b' - b_1^{2k+1} b'^{2k} = b_1 b'(b_2 - b_1^{2k} b'^{2k-1}) \geq 0 \rightarrow b_1 b' > 0$ thus $b_2 - b_1^{2k} b'^{2k-1} \geq 0$ (21)

Since based on lemma 1.5 in equation $a^{2k+1} + b^{2k+1} = c^{2k+1}$ the amount of $a+b-c$ is always positive and on the other hand the amounts of $a_1, b_1$ and $c_1$ are considered positive as well, so it can be concluded that in relations (19) to (21) the expression $c_1^{2k} c'^{2k-1} - c_2$, $a_2 - a_1^{2k} a'^{2k-1}$, and $b_2 - b_1^{2k} b'^{2k-1}$ are always positive.

**Theorem 1.7.** If the determinant of the system of three homogenous equations with three unknowns is equal to zero then each of the equations can be written in the form of a linear combination of two other equations.

Consider the system of three equations with three unknowns $a_1 x + b_1 y + c_1 z = 0$, $a_2 x + b_2 y + c_2 z = 0$, $a_3 x + b_3 y + c_3 z = 0$, assuming that its determinant is zero. By computing the determinant of system we will have:

$Det = a_1 (b_2 c_3 - b_3 c_2) - b_1 (a_2 c_3 - a_3 c_2) + c_1 (a_2 b_3 - a_3 b_2) = 0$ (22)

Now it will be indicated that each of the above equations can be written in the form of a linear combination of two other equations. Suppose that for example the equation $a_1 x + b_1 y + c_1 z = 0$ obtained from linear combination of two other equations. It's means we have:

$$a_1x+b_1y+c_1z=0 \quad a_1=m*a_2+n*a_3 \quad (23)$$

$$a_2x+b_2y+c_2z=0 \Rightarrow b_1=m*b_2+n*b_3 \quad (24)$$

$$a_3x+b_3y+c_3z=0 \quad c_1=m*c_2+n*c_3 \quad (25)$$

Now we replace the amounts of $a$, $b$, and $c$ obtained from the relations (23) to (25) in the relation (22), we will have:

$$Det=(m*a_2+n*a_3)(b_2c_3-b_3c_2)-(m*b_2+n*b_3)(a_2c_3-a_3c_2)+(m*c_2+n*c_3)(a_2b_3-a_3b_2)=$$

$$m.\underbrace{(a_2b_2c_3-a_2b_3c_2-b_2a_2c_3+b_2a_3c_2+c_2a_2b_3-c_2a_3b_2)}_{0}+n.\underbrace{(a_3b_2c_3-a_3b_3c_2-b_3a_2c_3+b_3a_3c_2+c_3a_2b_3-c_3a_3b_2)}_{0}=0$$

As it is observed the coefficient of m and n became equal to zero and it means that the amounts of the above matrix determinant is independed of the amounts m and n and therefore the primary assumption is correct and in fact each equation can be written in the form of a linear combinations of two other equations.

This subject can be shown for the other equations in similar ways.

**Lemma 1.8.** In equation $a^{2k}+b^{2k}=c^m$ that $gcd(a, b, c)=1$ and c is even, it will necessarily be $m=1$.

**Proof.** $(2m+1)^{2k}+(2n+1)^{2k}=(4p+1)+(4q+1)=2^1[2(p+q)+1]=2^1u=(2^rc_1)^m \rightarrow m=1$

**Lemma 1.9.** In equation $a^{2k}+b^{2k}=c^{2k}$ that $gcd(a, b, c)=1$ and $a$, $b$, and $c$ are coprime, can be written in form of couple-relations $c^k-b^k=A_1^{2k}$ and $c^k+b^k=A_1'^{2k}$ that $A_1$ and $A_1'$ are coprime and $a=A_1*A_1'$.

## 2. Proof of Fermat's Last Theorem

Consider the following equation:

$$a^x+b^x=c^x \quad (26)$$

In relation (26), $a$, $b$, and $c$ are all positive integer numbers and relatively prime in pairs. We going to prove that equation only in lieu of $x \leq 2$ there will be an answer. The above equation will be examined in two states separately when $x$ is odd and even.

### 2.1. Examining the equation $a^x+b^x=c^x$ when $x$ is odd

Consider the following equation:

$$a^{2k+1}+b^{2k+1}=c^{2k+1} \quad (27)$$

Based on Lemma 1.1 it is suffices to prove Fermat's Last Theorem for 4 and for every odd prime $p \geq 3$. Therefore assume $2k+1$ is an odd prime number. Now by using of the fundamental Theorem 1.3 and Lemmas 1.4 to 1.7, it will be proved that equation (27) will have an answer just in the lieu of $k=0$.

According to Theorem 1.3, if equation (27) set up, then will be set up three couple-relations (13)-(14), (15)-(16), and (17)-(18) simultaneously that in mentioned equations is $a=a_1*a_2*a'$, $b=b_1*b_2*b'$ and $c=c_1*c_2*c'$, and numbers $a_2, b_1, b_2, c_1, c_2, a', b'$ and $c'$ are relatively prime in pairs and are all positive integer number and we have two state: or $a'=b'=c'=1$ and or exactly one of them is equal to $2k+1$ and the other two exponents should be equal to 1 necessarily. In other words at least two of the numbers among common factors $a', b'$, and $c'$ should be equal to $1$. Without loss of generality, it is considered that $a' = b' = 1$ and therefore c' is equal to $1$ or equal to $2k+1$.

Now based on Theorem 1.3 and knowing that the relations (13), (15) and (17) are set up simultaneously and with respect to the relations and concepts mentioned above, it can be written that:

$$\begin{cases} a+b = c_1^{2k+1}c'^{2k} = a_1a_2 + b_1b_2 & \Rightarrow a_1(a_2) + b_1(b_2) - c_1(c_1^{2k}c'^{2k}) = 0 \quad (28) \\ c-a = b_1^{2k+1}b'^{2k} = b_1^{2k+1} = c_1c_2c' - a_1a_2 & \Rightarrow a_1(a_2) + b_1(b_1^{2k}) - c_1(c_2c') = 0 \quad (29) \\ c-b = a_1^{2k+1}c'^{2k} = a_1^{2k+1} = c_1c_2c' - b_1b_2 & \Rightarrow a_1(a_1^{2k}) + b_1(b_2) - c_1(c_2c') = 0 \quad (30) \end{cases}$$

Therefore according to Theorem 1.3, if the relation (27) is set up, if so three relations (28), (29), and (30) must be set up simultaneously. Now we try to study the answers of three equations with three unknowns mentioned above that are in terms of the unknown's $a_1$ and $b_1$ and $c_1$ and prove that the system of above equations will be set up just in lieu of $k = 0$.

If the determinant of the system of above equations will not be zero then the equations have a trivial solution $a_1 = b_1 = c_1 = 0$ that is subsequently concluded that $a = b = c = 0$ that there is a contradiction. Thus in order to have an answer for the above equations system the matrix determinant of coefficient must be equal to zero that if so the system will have infinite answerers. Through direct calculation of determinant of above equations ((28) to (30)) it can be found out that the determinant of coefficients is always equal to zero. If the determinant of coefficients of equations are computed there will be:

$$Determinant_{Eq.\ 25,26,27} = a_2\left[b_1^{2k}(-c_2c') + b_2c_2c'\right] - b_2\left[a_2(-c_2c') + a_1^{2k}c_2c'\right] - c_1^{2k}c'^{2k}\left[a_2b_2 - a_1^{2k}b_1^{2k}\right]$$
$$= a_2c_2c'\left(b_2 - b_1^{2k}\right) + b_2c_2c'\left(a_2 - a_1^{2k}\right) - c_1^{2k}c'^{2k}\left(a_2b_2 - a_1^{2k}b_1^{2k}\right) \quad (31)$$

Now we multiply and then divide the terms of the relation (31) to expression $(a_1b_1c_1)$ and we will have:

$$Det = \frac{1}{a_1b_1c_1}\left[a_1a_2c_1c_2c'\left(b_1b_2 - b_1^{2k+1}\right) + b_1b_2c_1c_2c'\left(a_1a_2 - a_1^{2k+1}\right) - c_1^{2k+1}c'^{2k}\left(a_1b_1a_2b_2 - a_1^{2k+1}b_1^{2k+1}\right)\right] \quad (32)$$

Base on Lemma 1.6 the amounts $a_1\left(a_2 - a_1^{2k}\right)$, $b_1\left(b_2 - b_1^{2k}\right)$, and $c_1c'\left(c_1^{2k}c'^{2k-1} - c_2\right)$ are equal to $(a+b-c)$ and also it is known that, $c = c_1c_2c'$, $a = a_1a_2$, and $b = b_1b_2$, therefore by replacing these amounts in relation (32) finally will have:

$$Det = \frac{1}{a_1b_1c_1}\{ac(a+b-c) + bc(a+b-c) - (a+b)[ab - (c-a)(c-b)]\} =$$
$$\{\frac{1}{a_1b_1c_1} * (a+b-c) * [ac+bc-ac-bc]\} = 0 \quad (33)$$

Since matrix determinant of the coefficients is equals to zero, therefore based on Theorem 1.7, if we consider each consider the equation (30), it is obtained from the linear combination of equations (28) and (29) that we represent that in this state just in the lieu of $k = 0$ the system of above equations has an answer. Because from the linear combination of the equations (28) and (29) in terms of desired coefficients $m$, $n$ and equalizing it with the equation (30) we will have:

$$a_1 : \begin{cases} m*a_2 + n*a_2 = a_1^{2k} \\ b_1 : \begin{cases} m*b_2 + n*b_1^{2k} = b_2 \\ c_1 : \begin{cases} m*(c_1c')^{2k} + n*c_2c' = c_2c' \end{cases} \end{cases}$$

$$(m+n).a_2 = a_1^{2k} \qquad a_2 \geq a_1^{2k} \text{ thus should } 1 \geq m+n \quad (34)$$
$$\Rightarrow (1-m)b_2 = n.b_1^{2k} \qquad \Rightarrow b_2 \geq b_1^{2k} \text{ thus should } n \geq 1-m \quad (35)$$
$$(1-n)c_2 = m.c_1^{2k}c'^{2k-1} \qquad c_1^{2k}c'^{2k-1} \geq c_2 \text{ thus should } 1-n \geq m \quad (36)$$

According to Lemma 1.6 if the equation $(a_1a_2)^{2k+1} + (b_1b_2)^{2k+1} = (c_1c_2c')^{2k+1}$ will be set up in the lieu of each positive a, b, or c, then the inequalities $a_2 \geq a_1^{2k}$, $b_2 \geq b_1^{2k}$, and $c_1^{2k}c'^{2k-1} \geq c_2$ will be always set up, so in order to set up the equalities (34) to (36) the inequalities $1 \geq m+n$, $n \geq 1-m$, and $1-n \geq m$ should be set up simultaneously that it is concluded the equation $m+n=1$ which is always the only common part of the above inequalities should be set up. Therefore the mentioned equations will be written in the form of following equations:

$$\begin{cases} a_2 = a_1^{2k} \\ n(b_2 - b_1^{2k}) = 0 \Rightarrow \text{Since } gcd(a_2, a_1^{2k}) = gcd(b_2, b_1^{2k}) = gcd(c_2, c_1^{2k}, c') = 1 \text{ therefore must } k = 0 \\ m(c_1^{2k}c'^{2k-1} - c_2) = 0 \end{cases}$$

So just in lieu of $k=0$ the system of three homogenous equations with three unknowns mentioned above will have an answer. By studding the other states of linear combination in above equations we come to similar conclusions. In general firstly it was indicated that if the equation (27) is set up then the relations (13), (15), and (17) (that are same relations (28) to (30)) should be set up simultaneously as well and since it was shown the determinant of coefficients matrix of the system of above equations is equal to zero, by using of properties of the linear combinations and by applying Lemmas 1.5 and 1.6 and so Theory 1.7, it was proved that the system of above equations will have an answer just in lieu of $k=0$ and subsequently it can be concluded that the equation $a^{2k+1} + b^{2k+1} = c^{2k+1}$ will have an answer just in lieu of $k=0$.

Generally based on Theorems 1.3 and 1.7 and Lemmas 1.4 to 1.6, it was proved that the equation $a^{2k+1} + b^{2k+1} = c^{2k+1}$ when $2k+1$ is prime number, only in recognition of $k=0$ will have an answer and based on Lemma 1.1 we can expand this statement for all odd numbers $2k+1$.

## 2.2. Examining the equation $a^x + b^x = c^x$ when x is even and equals to $4k$

Consider the following equation that $a$, $b$, and $c$ are integers and relatively prime in pairs:

$$a^{4k} + b^{4k} = c^{4k} \qquad (37)$$

Equation (37) can be written as follows:

$$(a^{2k} + b^{2k})^2 - (c^{2k})^2 = 2a^{2k}b^{2k} \qquad (38)$$

According to Lemma 1.8 in relation (37), $c$ is always odd, so $a$ must be odd and $b$ must be even or vice versa $a$ must be even and $b$ must be odd. Suppose that $a$ is odd and equals to $a = a_1.a_2$ and $b$ is even and equals to $b = 2^p * b_1 * b_2$ ($b_1$ and $b_2$ are odd and $p$ is indefinite). Thus both expressions $(a^{2x} + b^{2x} + c^{2x})$ and $(a^{2x} + b^{2x} - c^{2x})$ are even. Therefore relation (38) can be written in form of the relations (39):

$$(a^{2k} + b^{2k} - c^{2k})(a^{2k} + b^{2k} + c^{2k}) = 2(2^p * a_1 * a_2)^{2k}(b_1 * b_2)^{2k} = 2^{2kp+1} * a_1^{2k} * a_2^{2k} * b_1^{2k} * b_2^{2k} \qquad (39)$$

In relation (39), base on Lemma 1.8, expression $\left(a^{2k}+b^{2k}+c^{2k}\right)$ has factor $2^1$ (therefore expression $\left(a^{2k}+b^{2k}-c^{2k}\right)$ has factor $2^{2kp}$), so because $a$, $b$, and $c$ are coprime therefore base on Lemmas 1.8 and 1.9, the relation (39) can be written as follows:

$$a^{2k}+b^{2k}+c^{2k}=2a_1^{2k}*b_1^{2k} \tag{40}$$

$$a^{2k}+b^{2k}-c^{2k}=2^{2pk}a_2^{2k}*b_2^{2k} \tag{41}$$

That according to Lemma 1.9 in above relations parameters $a_1$, $b_1$, $a_2$, and $b_2$ must be relatively prime in pair. By adding the relations (40) and (41) we will have:

$$a^{2k}+b^{2k}=a_1^{2k}*b_1^{2k}+2^{2pk-1}a_2^{2k}*b_2^{2k} \tag{42}$$

As assume $a=a_1.a_2$ and $b=2^p.b_1.b_2$, so the relation (42) can be written as follows:

$$a_1^{2k}\left(a_2^{2k}-b_1^{2k}\right)=2^{2pk-1}b_2^{2k}\left(a_2^{2k}-2b_1^{2k}\right) \tag{43}$$

In the above equation as $a_1$ and $b_1$ are odd and coprime, thus if assume $t\,|\,\left(a_1^{2k}-2b_2^{2k}\right)$, then should $t\,|\,\left(a_2^{2k}-b_1^{2k}\right)$ and will have:

$$a_2^{2k}-b_1^{2k}=2^{2pk-1}b_2^{2k}t \tag{44}$$

$$a_2^{2k}-2b_1^{2k}=a_1^{2k}t \tag{45}$$

By solving two equations (44) and (45), the amounts $a_1^{2k}$ and $b_2^{2k}$ will be equal to:

$$b_1^{2k}=\left(2^{2kp-1}b_2^{2k}-a_1^{2k}\right)t \tag{46}$$

$$a_2^{2k}=\left(2^{2kp}b_2^{2k}-a_1^{2k}\right)t \tag{47}$$

As $gcd(a,b)=1$ and we suppose $b=b_1*b_2$ and $a=2^p*a_1*a_2$, so $gcd(b_1,a_2)=1$, then when just $t=1$ the condition of $gcd(b_1,a_2)=1$ will be set up, therefore the relations (46) and (47) will be formed as follows:

$$b_1^{2k}+a_1^{2k}=2^{2kp-1}b_2^{2k} \tag{48}$$

$$a_1^{2k}+a_2^{2k}=2^{2kp}b_2^{2k} \tag{49}$$

According to Lemma 1.8 in the relation (48) there must be $2kp-1=1\,(kp=1)$ and in relation (49) there must be $2kp=1\,(kp=0.5)$, firstly because p and k are both integer numbers, so $kp$ cannot be equal to $0.5$ and this is a contradiction and secondly as two relations mentioned above must be set up simultaneously, $kp$ cannot possess both amounts $0.5$ and $1$, therefore another contradiction is seen.

Generally it was proved that the equation $a^x+b^x=c^x$ when $x$ is even and equal to $4k$ will have no answer.

### 2.3. Examining the equation $a^x+b^x=c^x$ when $x$ is even and equals to $4k+2$

Consider the following equation:

$$a^{4k+2} + b^{4k+2} = c^{4k+2} \implies (a^2)^{2k+1} + (b^2)^{2k+1} = (c^2)^{2k+1} \tag{50}$$

As it was proved formerly in part (2-1), equation (50) will have an answer just in the lieu of $k=0$ that in this state, the above equation will be in the following form:

$$a^2 + b^2 = c^2 \tag{51}$$

Now we study the possible answers for the above equation. According to Lemma 1.8, $c$ must be always odd and among $a$ and $b$ one of them is even and the other one is odd. Suppose $a$ is even and $b$ is odd and equal to $b_1 * b_2$ ($b = b_1 * b_2$). In that case based on Lemmas 1.9, the equation (51) can be written in the form of the following couple-relations:

$$c - a = b_1^2 \tag{52}$$
$$c + a = b_2^2 \tag{53}$$

So the general answer of the equation (51) ($b$ is odd and equal to $b_1 * b_2$) will be as follows:

$$c = (b_1^2 + b_2^2)/2, \ a = (b_1^2 - b_2^2)/2, \ b = b_1 * b_2$$

Generally it was proved that the equation $a^{2k+1} + b^{2k+1} = c^{2k+1}$ only in lieu of $k=0$ will have an answer and the equation $a^{2k} + b^{2k} = c^{2k}$ just in lieu of $k=1$ has an answer. Thus it was confirmed that the equation $a^x + b^x = c^x$ just in lieu of $x \leq 2$ has an answer and this proves the correctness of Fermat's Last Theorem.

**Acknowledgment**

We would like to express our gratitude to Dr. Ali Rajaei, Dr. John Hamman, and Dr. Reza Taleb, for providing their support and encouragement.

**3. Conclusion**

Generally for the equation $a^x + b^x + c^x = 0$ by using properties of the algebra identities and linear algebra and presenting a new theorem, it was represented that: If $x$ is odd and equal to $2k+1$, it was proved according to key Theorem 1.3 and Lemma 1.7 that just in lieu of $k=0$ the equation $a^{2k+1} + b^{2k+1} = c^{2k+1}$ will have an answer. So if $x$ is even and equal to $2(2k)$ if so the equation $a^{4k} + b^{4k} = c^{4k}$ will have no answer and if $x$ is even and equal to $2(2k+1)$, based on the first result in above, it can be concluded that the equation $(a^2)^{2k+1} + (b^2)^{2k+1} = (c^2)^{2k+1}$ will have an answer just in lieu of $k=0$ and its general answer for odd p and odd q will be: $c = (p^2 + q^2)/2, \ a = (p^2 - q^2)/2, \ b = p*q$. Based on the above conclusions it was proved that the equation $a^x + b^x = c^x$ just in lieu of $x \leq 2$ will have an answer.